\documentclass[a4paper,12pt]{amsart}
\usepackage[utf8x]{inputenc}
\usepackage{amsmath,amssymb,amsthm}
\usepackage{color}
\definecolor{myblue}{rgb}{0.09,0.32,0.44} 
\usepackage
	[
	colorlinks,
        linkcolor=myblue,
        citecolor=myblue,
        urlcolor=blue,
        pdfborder={0 0 0},
        pdfborderstyle={},
	]
{hyperref}
\usepackage[all]{xy}
\usepackage{float}
\usepackage{breakurl}

\newcommand{\ZZ}{\mathbb{Z}}

\newcommand{\Ht}{{\rm Ht}}
\newcommand{\red}{{\rm red}}

\newtheorem*{theorem*}{Theorem}

\begin{document}

\title[Irreducible $\pm1$ polynomials]{Is a bivariate polynomial with $\pm 1$ coefficients irreducible? Very likely!}
\author{Lior Bary-Soroker}
\address{School of Mathematical Sciences, Tel Aviv University, Ramat Aviv, Tel Aviv 6997801, Israel. e-mail: \tt{barylior@post.tau.ac.il}}
\author{Gady Kozma}
\address{
Weizmann Institute, Rehovot, 76100, Israel. \newline
e-mail: \tt{gady.kozma@weizmann.ac.il}
}
\maketitle

\section{Introduction}
The \emph{height} of a polynomial $f(X) = \sum_{i=0}^r a_i X^i$ with integral coefficients is defined as 
\[
\Ht(f) = \max_{i} |a_i|.
\]
The probability that a polynomial $f\in \ZZ[X]$ of degree $r$ and height at most $N$ to be reducible is 
\begin{equation}\label{eq:prob_red}
\frac{1}{2N+1}\leq p_{r,N} \leq C_r \frac{1}{N},
\end{equation}
for some constant $C_r>0$ depending only on $r$. 
The left hand side follows from the fact that if $a_0=0$; i.e., $f(0)=0$, then the polynomial is reducible. Upper bounds as in the left hand side were studied for a long period of time, see e.g.\ \cite{Dorge,Gallagher,Kuba,Rivin}. In particular, as $r$ is fixed and $N$ tends to infinity, we get the rate of decay of $p_{r,N}$.

One would expect that as $N$ being fixed and as $r$ tends to infinity and unless there is some obvious obstruction, we would still get that 
\[
\lim_{r\to \infty} p_{r,N}=0.
\] 
More generally, take a finite set $S$; say $S=\{\pm 1\}$, and let
$p_{r,S}$ be the probability that a random polynomial of degree $r$
with coefficients in $S$ is reducible. It is notoriously difficult to prove that 
\[
\lim_{r\to \infty} p_{r,S} =0. 
\] 
Some cases that were studied extensively in the literature are
$S=\{0,1\}$ \cite{konyagin} and $S=\{\pm 1\}$ (Littlewood polynomials)
\cite{mathoverflow}. See also the realted \cite{PSZ}. We have no solution for this problem. 

Instead, the goal of this study is much more modest. We address a much simpler question when we add one degree of freedom; namely, we replace univariate polynomials by bivariate polynomials. We hope that people working in this area might be interested in this case. We prove:

\begin{theorem*}\label{thm:main}
Let $F=F(X,Y) = \sum_{i,j\leq r} \pm X^iY^j$ be a bivariate polynomial of degree $r$ with random coefficients $\pm 1$. 
Then 
\[
\lim_{r\to \infty}\mathbb{P}(F \mbox{ \rm reducible}) = 0.
\]
\end{theorem*}

The key tool in the proof of the theorem is a variant of
\eqref{eq:prob_red} for polynomials with odd coefficients with an
extra $(\log N)^2$ factor and with an explicit bound in term of $r$: Put 
\begin{equation}\label{eq:defOme}
\Omega(r,N) = \left\{f = \sum_{i=0}^r a_i X^i :   a_r\equiv 1\pmod 2 , \Ht(f)\leq 2N-1\right\}.
\end{equation}
Then there exists an absolute constant $C>0$ such that for any $r>1$
and $N>2$ the probability that a random uniform polynomial $f\in \Omega(r,N)$ is reducible satisfies
\begin{equation}\label{eq:odd}
\mathbb{P}_{\Omega(r,N)}( f \mbox{ reducible}) \leq C\cdot \frac{r(\log N)^2}{N}
\Big(1+\frac 1{2N}\Big)^r.
\end{equation}

\section{Proof of \eqref{eq:odd} -- Rivin's argument}
Our proof of \eqref{eq:odd} is based on the surprising idea of Rivin \cite{Rivin} who noticed that for monic polynomials conditioning on the free coefficient gives a saving of $\log N/N$. In our case the polynomials are non-monic; hence we condition also on the leading coefficients, whence the $(\log N)^2$ factor:

Put $\Omega=\Omega(r,N)$ (from (\ref{eq:defOme})) and note that as there are $2N$ odd integers in $[-2N+1, 2N-1]$,  one has
\begin{equation}\label{eq:Omega}
\#\Omega= (2N)^{r+1}.
\end{equation}  
Let us fix $s,t>0$ with $r=s+t$ and $b_0,c_0,b_s,c_t\in \ZZ$ such that $a_0=b_0c_0$  and $a_{r} = b_sc_t$ are odd integers with $|a_0|,|a_r|\leq 2N-1$. Let $V=V(b_0,c_0,b_s,c_t,s,t)$ be the set of $f\in \Omega$ such that $f=gh$ with $\deg g=s$, $\deg h=t$, $g(0)=b_0$, $h(0)=c_0$, and the leading coefficients of $g,h$ are $b_s,c_t$, respectively. 

Put $M=2N+1$ and consider the reduction of coefficient map  $\red\colon \Omega \to \mathbb{Z}/M\mathbb{Z}[X]$. As $M>2N$ and $\gcd(M,2)=1$, it follows that $\red$ is injective. Moreover, the image $\red(V)$ of $V$ is contained of  in the set of  products $\bar{g}\bar{h}$, where $\bar{g}$ (resp., $\bar{h}$) ranges over polynomials of degree $s$ (resp., $t$) with coefficients in $\ZZ/M\ZZ$, leading coefficients $b_s\mod M$ (resp., $c_t\mod M$), and such that $\bar{g}(0) \equiv b_0 \mod M$ (resp., $\bar{h}(0) \equiv c_0\mod M$). Since there are $M^{s-1}M^{t-1} = M^{r-2}$ pairs $(\bar{g}, \bar{h})$ of such polynomials, we get 
\begin{equation}\label{eq:a_0neq0}
\#V(b_0,c_0,b_s,c_t,s,t) \leq M^{r-2}.
\end{equation}
Each reducible polynomial $f\in \Omega$ with $a_0=f(0)$ and leading coefficient $a_r$ lies in $V(b_0,c_0,b_s,c_t,s,t)$ for some parameters satisfying $b_0c_0=a_0$, $b_sb_t = a_r$, $s,t\geq 1$, and $s+t=r$. Thus, by 
\eqref{eq:a_0neq0}, one has
\begin{equation}\label{sum}
\begin{split}
\#\{ f\in &\Omega(r,N) : f \mbox{ reducible}\} \\
& \qquad \leq 
\sum_{a_0,a_r} \sum_{b_0\mid a_0,\ b_s\mid a_r} \sum_{s+t = r} \# V(b_0,\frac{a_0}{b_0}, b_s,\frac{a_r}{b_s},s,t) 
\\ & \qquad\leq (r-1)M^{r-2}
\sum_{b_0,b_s}2\left\lfloor\frac{2N-1}{b_0}\right\rfloor
\cdot 2\left\lfloor\frac{2N-1}{b_s}\right\rfloor.
\end{split}
\end{equation}
Here $a_0,a_r$ run over odd integers with $|a_0|, |a_r|\leq 2N-1$. 
A simple calculation bounds the last sum by $CM^2(\log M)^2$. Together with \eqref{sum} and \eqref{eq:Omega} one directly gets  \eqref{eq:odd} (the difference between $M^{r}$ and $(2N)^{r}$ is responsible for the last term in \eqref{eq:odd}). \qed

\section{Proof of the Theorem}
Notice that if $F(X,Y)$ is reducible, then either $F(2,Y)$ is
reducible, or $F(X,2)$ is reducible, or $F(X,Y) = f(X)g(Y)$. 
Indeed, assume $F(X,Y)= f(X,Y)g(X,Y)$, let $r=\deg_Y F$, $s=\deg_Y f$, and $t=\deg_Y g$. Then $r=s+t$. Substituting $X=2$, gives $F(2,Y) = f(2,Y)g(2,Y)$. The coefficient of $Y^r$ is $\pm1\pm 2\pm  \cdots \pm 2^r$, hence odd and in particular nonzero. Thus $\deg F(2,Y)=r$. Since $\deg f(2,Y)\le s$ and $\deg g(2,Y)\le t$ we may write
\[
r=\deg F(2,Y)=\deg f(2,Y)+\deg g(2,Y)\le s+t=r
\]
so both inequalities are in fact equalities, i.e.\ $\deg f(2,Y)=s$ and $\deg g(2,Y)=t$. 
 Thus, either $F(2,Y)$ is reducible, or one of the pair $s,t$ is zero; say $s=0$. Then $f(X,Y)=f(X)$. Similarly, by substituting $Y=2$, we get that either $F(X,2)$ is irreducible or $g(X,Y)=g(Y)$, as claimed.

As
\[
F(2,Y) = \sum_{i=0}^r (\sum_{j=0}^{r} \pm 2^j) Y^i,
\]
the coefficient of $Y^i$ is a random odd number in $[-2^{r+1} +1, 2^{r+1}-1]$.
%We note that if $F(X,Y)$ is reducible, then so is $F(2,Y)$. 
Thus, by \eqref{eq:odd} with $N=2^r$, we get
\[
\mathbb{P}(F(2,Y) \mbox{ reducible}) \ll \frac{r^3}{2^r} .
\]
($A\ll B$ is short for $A=O(B)$). Similarly, $\mathbb{P}(F(X,2) \mbox{ reducible}) \ll \frac{r^3}{2^r}$. If $F(X,Y)= f(X)g(Y)$, then both $f$ and $g$ can be normalized to have $\pm1$ coefficients; hence 
\[
\mathbb{P}(\exists f,g\textrm{ s.t. }F(X,Y) = f(X)g(Y)) \leq \frac{2^{r+1}\cdot 2^{r+1}}{2^{(r+1)^2}}\ll 2^{-r^2}.
\]
We conclude that
\[
\mathbb{P}(F(X,Y) \mbox{ reducible}) \ll \frac{r^3}{2^r} \to 0, \qquad r\to \infty,
\]
as needed.
\qed

\subsection*{Acknowledgements}
We thank Benjamin L.\ Weiss, Josh Zelinsky and the anonymous referee for pointing out gaps in earlier versions and Mark Shusterman for helpful discussions.

LBS is partially supported by the Israel Science Foundation grant no.~40/14.
GK is partially supported by the Israel Science Foundation grant
no.~1369/15 and by the Jesselson Foundation.


\begin{thebibliography}{A}
\bibitem{Dorge}
K. D\"orge,
\newblock Absch\"atzung der Anzahl der reduziblen Polynome,
\newblock Math. Ann. 160 (1965), 59--63.
\newblock Available at:
\href{http://link.springer.com/article/10.1007%2FBF01364334}{\nolinkurl{springer.com/10.1007}}


\bibitem{Gallagher}
Patrick X. Gallagher,
\newblock The large sieve and probabilistic Galois theory, 
\newblock in: Proc. Sympos. Pure Math., Vol. XXIV, Amer. Math. Soc. (1973), 91--101.
\newblock Available at: \burl{http://bookstore.ams.org/pspum-24/}

\bibitem{konyagin} 
Sergei V. Konyagin, 
\newblock On the number of irreducible polynomials with 0,1 coefficients.
\newblock Acta Arith. 88:4 (1999), 333--350. 
\newblock Available at:
\href{http://pldml.icm.edu.pl/pldml/element/bwmeta1.element.bwnjournal-article-aav88i4p333bwm?q=bwmeta1.element.bwnjournal-number-aa-1999-88-4;3&qt=CHILDREN-STATELESS}{\nolinkurl{pldml.icm.edu.pl/aa8844}}

\bibitem{Kuba}
Gerald Kuba, 
\newblock On the distribution of reducible polynomials,
\newblock Math. Slovaca 59 (2009), no.\ 3, 349--356.
\newblock Available at:
\href{http://link.springer.com/article/10.2478%2Fs12175-009-0131-6}{\nolinkurl{springer.com/10.2478}}

\bibitem{PSZ}
Ron Peled, Arnad Sen and Ofer Zeitouni,
\newblock Double roots of random Littlewood polynomials,
\newblock Israel J. Math, in press.
\newblock Available at: \href{http://arxiv.org/abs/1409.2034}{\nolinkurl{arXiv:1409.2034}}




\bibitem{Rivin}
Igor Rivin, 
\newblock  Galois Groups of Generic Polynomials,
\newblock Preprint.
\newblock \href{http://arxiv.org/abs/1511.06446}{\nolinkurl{arXiv:1511.06446}}

\bibitem{mathoverflow} 
Some guy on the street, 
\newblock Irreducible polynomials with constrained coefficients,
\newblock MathOverflow.
\newblock Available at: \href{http://mathoverflow.net/questions/7969/\\irreducible-polynomials-with-constrained-coefficients}{\nolinkurl{mathoverflow.net/7969}}
\end{thebibliography}
\end{document}